\documentclass[11pt,reqno,fleqn,a4paper]{article}
\usepackage{amsmath,amssymb,amsfonts,amsthm,enumitem,cite}
\usepackage{setspace,xcolor}
\usepackage{graphicx}
\usepackage{multicol}
\usepackage{lipsum,lineno}
\usepackage{ragged2e}
\usepackage{multirow}
\usepackage{bm}
\usepackage[margin=1.0in]{geometry}
\usepackage[utf8]{inputenc}
\usepackage[T1]{fontenc}
\usepackage{mathtools}
\usepackage[font=small,labelfont=bf]{caption}
\usepackage[subrefformat=parens]{subcaption}
\usepackage{mwe}
\usepackage{mfirstuc}

\providecommand{\keywords}[1]{\textbf{\textbf{\\Keywords. }}#1}
\providecommand{\class}[1]{\textbf{\textbf{Mathematics Subject Classification (2010). }}#1}
\usepackage{cite}
\MFUnocap{are}
\MFUnocap{or}
\MFUnocap{etc}

\def\capmystringaux#1#2\relax{\uppercase{#1}\lowercase{#2}}
\usepackage{lineno,hyperref}
\newcommand{\subfigANDtitle}[2][.2\linewidth]{%
	\begin{tabular}{@{}>{\centering\arraybackslash}p{#1}@{}} #2 \end{tabular}}
\title{}
\author{}
\DeclareSymbolFont{slenderlargesymbols}{OMX}{ccex}{m}{n}
\usepackage{amsmath,amssymb,amsfonts,calrsfs,bbm}
\usepackage{setspace}
\usepackage{graphicx,array}
\usepackage{amsthm}
\usepackage{multicol}
\usepackage{hyperref}
\usepackage{lipsum,lineno}
\usepackage{amsthm}
\usepackage{ragged2e}
\usepackage{multirow}
\usepackage{bm}
\usepackage{mathtools}
\usepackage{mwe}
\usepackage{breqn}
\usepackage{mfirstuc}
\MFUnocap{are}
\MFUnocap{or}
\MFUnocap{etc}
\usepackage{tikz,caption,subcaption}
\usetikzlibrary{shapes}
\usetikzlibrary{positioning}
\tikzset{main node/.style={circle,fill=black,draw,minimum size=.2cm,inner sep=0pt},
}
\tikzset{node1/.style={diamond,draw,minimum size=.2cm,inner sep=0pt},
}
\tikzset{node2/.style={star,draw,minimum size=.2cm,inner sep=0pt},
}

\DeclarePairedDelimiterX{\Set}[2]\{\}{%
	\, #1 \;\delimsize\vert\; #2 \,
}

\usepackage{amsfonts}
\usepackage{mathrsfs}

\newtheorem{theorem}{Theorem}[section]

\newtheorem{proposition}[theorem]{Proposition}
\theoremstyle{definition}
\newtheorem{definition}[theorem]{Definition}
\theoremstyle{remark}
\newtheorem{remark}[theorem]{Remark}
\newtheorem*{example}{Example}
\numberwithin{equation}{section}
\date{}
\begin{document}
	\title{\textbf{ On the convergence of sequences in the space of $n$-iterated function systems with applications}}
	\author{Praveen M\footnote{E-mail:praveenzmedium@gmail.com, praveenm@nitc.ac.in}, Sunil Mathew\footnote{E-mail:sm@nitc.ac.in}\\ \small Department of Mathematics,\\\small National Institute of Technology Calicut, Calicut - 673 601\\\small India.
}
	\maketitle	
\begin{abstract}
	
This article discusses the notion of convergence of sequences of iterated function systems. The technique of iterated function systems is one of the several methods to construct objects with fractal nature, and the fractals obtained with this method are mostly self-similar. The progress in the theory of fractals has found potential applications in the fields of physical science, computer science, and economics in abundance. This paper considers the metric space of $n$- iterated function systems by introducing a metric function on the set of all iterated function systems on a complete metric space consisting of $n$ contraction functions.
Further, sequences of $n$- iterated function systems with decreasing, eventually decreasing, Cauchy and convergent properties are discussed. Some results on sequences of $n$- iterated function systems and sequences of contractions are obtained. The practical usage of the theory discussed in the article is explored towards the end.
	
\medskip \medskip
	
	\noindent\class{Primary 28A80; Secondary
		11B05.}\medskip
\end{abstract}

\section{Introduction}
The fractal theory is a leading research area of mathematics which has made its own identity in many of the interdisciplinary sciences. Recently research conducted by MIT scientists in condensed matter physics discovered\cite{li2019scale} fractal patterns in neodymium nickel oxide($NdNiO_3$), a quantum material that is rare in the earth. The quantum, atomic-scale effects of the quantum materials result in the bizarre electronic or magnetic behavior of them. The quantum material $NdNiO_3$, depending on its temperature, behaves both as an electrical conductor and an insulator. The analysis conducted by the researchers with the aid of statistics of domain distribution on the texture of the magnetic domains of $NdNiO_3$ helped them identify a fractal pattern in it. The scientific community is conducting an extensive study on $NdNiO_3$ for the immense applications it offers, such as the possibility to use it as a building block for neuromorphic devices, which are the artificial systems that imitate biological neurons. The researchers believe the knowledge regarding the nanoscale magnetic and electronic textures is central to study and engineer other materials for similar scopes.

Fractal geometry, which was established in the mid-1980s with the help of computers, over the years, became a bridge in the gap between pure mathematics and applied sciences. 
Benoit B. Mandelbrot, who is considered the founding father of this branch of study, laid down the basic structure for fractal theory and made many remarkable contributions in this field.
 Mandelbrot brought together concepts in the analysis by people like Felix Hausdorff, Pierre Fatou, Gaston Julia, and concepts in geometry by people like Helge von Koch, Ernesto Cesaro, in the construction of the theory of fractals. Self-similarity, which was defined by Cesaro in 1905, Hausdorff dimension, defined by Hausdorff in 1918, and non-differentiability are some of the crucial notions used in building the theory of fractals by Mandelbrot. He used the term `fractal' to describe repeating or self-similar mathematical patterns. He also gave a formal definition of fractals in 1975, comparing topological dimension and Hausdorff dimension of an object, which later he retracted, claiming that there existed objects which should be considered as fractals and not qualify as a fractal with his definition. The mathematicians are still working on finding a formal definition of fractals. In his celebrated book Fractal Geometry of Nature(1982), Mandelbrot highlighted the many occurrences of objects in nature with the fractal properties.
 Mandelbrot claimed that the notions from fractal theory could be used to comprehend the components of any essential structure in nature and make predictions about their future.
 
There are many natural objects as well as phenomena following the fractal characteristics. The shape of coastlines, the venation of leaves and branching of trees, the branching of blood vessels and nerves in the human body, the DNA molecule, the price history in the stock market are some of the places where we see fractal properties. The fractal theory has many real-life applications, as well. The traders often make use of fractals to understand the direction in which the price will develop in the stock market. The new researches conducted in technology established that antennae having specific fractal shapes will reduce, on a considerable scale, the size and weight of the antennae. The fractal patterns derived from our blood vessels have been used in the silicon chips of computers to allow the cooling fluid to uniformly flow across the surface of the chip and keep it cool. Many branches of applied sciences are making use of the fractal theory in their research to revolutionize life in the future.

 There are several methods to generate a fractal object mathematically. The widely used methods include iterated function systems(IFS), strange attractors, L-systems, escape-time fractals, random fractals, and finite subdivision rules.
%
%
The technique of iterated function systems generate fractal objects with self-similarity; that is, the part of the object resembles the whole.
 The concept of iterated function system(IFS) first appeared in a paper\cite{hutchinson1979fractals} by John E. Hutchinson in 1981 and later was popularized by Michael Barnsley in his book\cite{barnsley2014fractals}. According to Hutchinson IFS theory, an iterated function system consists of a family of Banach contracting self-maps on a complete metric space. He proved that the Hutchinson operator or the fractal operator defined by him on the hyperspace of all non-empty compact subsets of the underlying space with the Hausdorff-Pompieu metric exhibits a unique fixed point namely the attractor of the corresponding IFS.
%
 Following this, mathematicians have extended the Hutchinson IFS theory by using more general spaces as the underlying space, taking infinite number of contractions in the IFS, and using generalized contraction mappings instead of Banach contractions. The generalization in  which an infinite set of contractions is used in place of a finite set is called an \textit{infinite iterated function system}($IIFS$)\cite{fernau1994infinite} and is called a \textit{countable iterated function system}($CIFS$)\cite{secelean2013countable,secelean2012existence} when the set of contractions is countable. Also, there are generalizations in the literature with the contraction condition relaxed to generalized contraction conditions such as $r$-contraction\cite{secelean2013iterated}, convex contraction\cite{istratescu1981some,istractescu1982some,miculescu2015generalization,georgescu2017ifss}, Meir-Keeler type contraction, $F$-contraction\cite{secelean2013iterated}, weak contraction\cite{hata1985structure}, etc. The other types of generalizations include the relaxation on the completeness condition of the underlying metric space, replacing the metric space with a product of metric spaces, relaxing the metric condition in the space to partial metric\cite{minirani2014fractals} condition, 
 replacing the metric space with a general topological space\cite{mihail2012topological}, etc.
 
In this paper, we introduce the space of iterated function systems consisting of a fixed number of contractions. We also give an ordering of an IFS with respect to another IFS in this space. Further, it is defining several types of sequences of IFSs, such as decreasing, eventually decreasing, Cauchy, convergent, etc. Using this basic structure, we obtain certain results connecting the attractors of IFSs in a sequence of IFS.
A possible application of the theory developed in the paper is discussed towards the end. 
\section{Preliminaries}
This section is providing the fundamental definitions and results that are the building blocks of the theory discussed in this paper, and they are from \cite{barnsley2014fractals,hutchinson1979fractals,hata1985structure}.
Throughout this paper, $(X, d)$ denote a complete metric space,  $\mathbbm{H}(X)$ denote the non-empty closed and bounded subsets of $(X, d)$ and $\mathbbm{K}(X)$ denote the non-empty compact subsets of $(X, d)$.
  

 
The following are the concepts from the literature that help in defining a metric called the Hausdorff metric in $\mathbbm{H}(X)$.
\begin{definition}\cite{barnsley2014fractals}
	Let $(X,d)$ be a metric space, $x\in X$ and $K,A,B\in \mathbbm{H}(X)$. Then $$d(x,K)=\inf\{d(x,y) : y\in K\} \text{ and}$$
	$$d(A,B)=\sup\{d(x,B) : x\in A\}.$$
\end{definition}
The function $d:\mathbbm{H}(X)\times \mathbbm{H}(X) \to \mathbb{R} $ is not a metric because of the the following remark.
\begin{remark}\cite{barnsley2014fractals}
	In general, $d(A,B)\ne d(B,A)$ and $d(A,B)=0$ even if $A\ne B$.
\end{remark}
The Hausdorff metric on $X$ is used to measure how far any given two subsets of $X$ are. It is the largest of all the distances from a point in one set to the nearest point in the other set. The definition of Hausdorff metric is as follows.
\begin{definition}\cite{barnsley2014fractals} The map $h:\mathbbm{H}(X)\times \mathbbm{H}(X)\to \mathbb{R} $ defined by
	$$h(A, B) = \max \{d(A,B), d(B,A)\}$$ is a metric on $\mathbbm{H}(X)$ called the\textit{ Hausdorff-Pompeiu metric} or \textit{ Hausdorff metric}. The metric space $(\mathbbm{H}(X),h)$ is complete provided that $(X,d)$ is complete.
\end{definition}
The definition of a dynamical system is as follows:
\begin{definition}\cite{barnsley2014fractals}
	A dynamical system is a transformation $f:X\to X$ on a metric space $(X,d)$. It is denoted by $\{X;f\}$.
\end{definition}
An iterated function system is a dynamical system with finitely many contraction maps acting on a complete metric space. The following is the formal definition of a hyperbolic $IFS$.
\begin{definition}\cite{barnsley2014fractals}
	A \textit{hyperbolic Iterated Function System }(hyperbolic IFS) consists of a complete metric space $(X,d)$ and a finite number of contraction mappings $f_i:X\to X$, with respective contractivity factors $t_i$ for $i=1,2,\cdots,n.$ The value $t=\max_{i=1}^nt_i$ is called the contractivity factor of the IFS.
	
\end{definition}
The $IFS$ is a method to construct objects with fractal nature mathematically. The attractor or the set fixed point obtained from the $IFSs$ are mostly self-similar fractals. We have the following theorem for the existence and uniqueness of the attractor of an $IFS$.
\begin{theorem}\cite{barnsley2014fractals}
	Let $\displaystyle \Bigl\{X;f_i,i=1,2,\cdots,n\Bigr\}$ be an IFS with contractivity factor $t$. Then the transformation $W:\mathbbm{H}(X) \to \mathbbm{H}(X)$ defined by $\displaystyle W(B)=\cup_{i=1}^nf_i(B)$ for all $B\in \mathbbm{H}(X)$, is a contraction mapping on the complete metric space $(\mathbbm{H}(X),h)$ with contractivity factor $t$. Its unique fixed point, $A\in \mathbbm{H}(X)$, exists and is given by $\displaystyle A=\lim_{n\to \infty}W^{[n]}(B)$ for any $B\in \mathbbm{H}(X)$.
\end{theorem}
To get a mathematical model for a given object with fractal nature one can make use of $IFS$ theory. The following theorem called the \textit{Collage theorem} ensures how this can be achieved.
\begin{theorem}\cite{barnsley2014fractals}
	Let $L\in \mathbbm{H}(X)$ and $\epsilon>0$ be given. Choose an $IFS$ $\displaystyle\Big\{X;f_i,i=1,2,\cdots,n\Bigr\}$ with contractivity factor $0\le t<1$, so that $h\Big(L,\displaystyle\bigcup_{i=1}^{n}f_i(L)\Big)\le \epsilon$
	Then $h\Big(L,F\displaystyle\Big)\le \frac{\epsilon}{1-t},$
	where $F$ is the attractor of the IFS.
\end{theorem}
In the next section we provide our definitions and major results.
\section{Definitions and Main Results}
This section provides new definitions, associated results and the impact of them on the existing theory. We first make the set of all iterated function systems on a complete metric space into a metric space by introducing a metric on it. In order to achieve this we need to structure the IFSs according to the given IFSs and the following definition explains how to do it.
\begin{definition}
	Consider a complete metric space $(X,d)$. Let $\mathcal{S}_X^n$ denotes the \linebreak collection of all iterated function systems on $X$ with $n$-contractions.
	Let \linebreak$\mathbb{S}=\{X;f_1,f_2,\cdots,f_n\}$ and $\mathbb{T}=\{X;g_1,g_2,\cdots,g_n\}$ be iterated function systems, i.e., $\mathbb{S},\mathbb{T}\in\mathcal{S}_X^n$. Then, we say, $\mathbb{T}$ \textit{is minimally ordered with respect to} $\mathbb{S}$, if $$\sum_{i=1}^n\bar{d}_\infty(f_i,g_i)=\min_{\sigma\in S_n}\sum_{i=1}^n\bar{d}_\infty(f_i,g_{\sigma(i)})$$
	where $S_n$ is the permutation group on $\{1,2,\cdots,n\}$ and $\displaystyle\bar{d}_{\infty}(f,g)=\sup_{x\in X}\frac{d(f(x),g(x))}{1+d(f(x),g(x))}$ is a metric on the set of all functions from $X$ to itself.
	
	Let $\mathbb{S},\mathbb{T}\in\mathcal{S}_X^n$, and $\mathbb{T}$ be minimally ordered with respect to $\mathbb{S}$. Then,
	we say, $\mathbb{S}\le\mathbb{T}$ if $c_{f_i}\le c_{g_i}$ for every $i=1,2,\cdots,n$, where $c_{f_i}, c_{g_i}$ denote the contractivity factors of $f_i,g_i$ respectively.
	
	Let $\mathbb{S}_j=\{X;f_{1_j},f_{2_j},\cdots, f_{n_j}\}$ be iterated function systems for $j=1,2,\cdots$ Then $(\mathbb{S}_j)_{j\ge1}$ is called a sequence of iterated function systems on $X$ with $n$-contractions. Let $c_{i_j}$ be the contractivity factor of the contraction $f_{i_j}$ for $1\le i\le n, j\ge 1$. The sequence of IFSs $(\mathbb{S}_j)_{j\ge1}$ is said to be decreasing if $\mathbb{S}_{j+1}\le \mathbb{S}_j$ for all $j=1,2,\cdots$ Also $(\mathbb{S}_j)_{j\ge1}$ is said to be eventually decreasing if there exists $k\in\mathbb{N}$ such that $\mathbb{S}_{j+1}\le \mathbb{S}_j$ for every $j\ge k$.
	
	The sequence $(\mathbb{S}_j)_{j\ge1}$ is a \textit{minimally ordered sequence of IFSs} if $\mathbb{S}_{j+1}$ is minimally ordered with respect to $\mathbb{S}_j$, for every $j=1,2,\cdots$
\end{definition}
\begin{remark}
	Hereafter an iterated function system on $X$ with $n$-contractions will be called an $n$-iterated function system on $X$.
	
\end{remark}
Now we define a decreasing and eventually decreasing sequence of contractions as below.
\begin{definition}
	Let $Con(X)$ denote the collection of all contractions on $X$. Then $\bar{d}_\infty$ is a metric on $Con(X)$. A sequence $(f_n)$ in $(Con(X),\bar{d}_\infty)$ is decreasing if $c_{n+1}\le c_n$, for every $n\ge 1$, where $c_n$ is the contractivity factor of $f_n$. Also, $(f_n)$ is eventually decreasing if there exists $N\in \mathbb{N}$ such that $c_{n+1}\le c_n$, for every $n\ge N$.
\end{definition}
\begin{remark}\label{con_rem}
	Even though $(X,d)$ is a complete metric space, $(Con(X),\bar{d}_\infty)$ need not be complete. 
\end{remark}
The following is an example which demonstrate the remark \ref{con_rem}.
\begin{example}
	Consider the metric space $X=[0,1]$ with Euclidean metric. Let $(f_n)_{n\ge1}$ be a sequence of contractions on $X$ defined by $f_n(x)=(1-\frac{1}{n})x$ and let $f(x)=x$. We have $$\begin{aligned}
	\bar{d}_\infty(f_n,f)&=\sup_{x\in X}\frac{d(f_n(x),f(x))}{1+d(f_n(x),f(x))}\\&=\sup_{x\in[0,1]}\frac{|(1-\frac{1}{n})x-x|}{1+|(1-\frac{1}{n})x-x|}\\
	&=\frac{\frac{1}{n}}{1+\frac{1}{n}}=\frac{1}{n+1}\to 0 \text{ as } n\to \infty.
	\end{aligned}$$
	Hence $f_n\to f$ in $\bar{d}_\infty.$ But $f$ is not a contraction. Therefore, $(Con(X),\bar{d}_\infty)$ is not complete.
\end{example}
The next theorem discusses the convergence of a sequence of contractions.
\begin{theorem}\label{con}
	Consider a complete metric space $(X,d)$. Let $(f_n)_{n\ge1}$ be an eventually decreasing Cauchy sequence in $(Con(X),\bar{d}_\infty)$. Then $(f_n)_{n\ge1}$ converges in $Con(X)$.
	\begin{proof}
		Since $(f_n)_{n\ge1}$ is eventually decreasing, we get $(c_n)_{n\ge 1}$ is eventually decreasing in $[0,1)$, where $c_n$ denotes the contractivity factor of $f_n$. Further we get $(c_n)_{n\ge1}$ converges in $[0,1)$, say, to $c$, because $(c_n)_{n\ge 1}$ is also bounded below by $0$, i.e., $\displaystyle\lim_{n\to \infty}c_n=c$.
		\\
		Also, since $(f_n)_{n\ge1}$ is Cauchy, we get, for any given $0<\epsilon<1$, there exists $N\in \mathbb{N}$ such that $\bar{d}_\infty(f_n,f_m)<\epsilon$ for every $n,m\ge N$, i.e.,
		$$\begin{aligned}
		&\sup_{x\in X}\frac{d(f_n(x),f_m(x))}{1+d(f_n(x),f_m(x))}<\epsilon, \forall n,m\ge N\\
		\implies& \frac{d(f_n(x),f_m(x))}{1+d(f_n(x),f_m(x))}<\epsilon, \forall x\in X \text{ and }n,m\ge N\\
		\implies&d(f_n(x),f_m(x))<\frac{\epsilon}{1-\epsilon},\forall x\in X \text{ and }n,m\ge N
		\end{aligned}$$
		Therefore, $f_n(x)$ is Cauchy in $(X,d)$, and $f_n(x)$ converges in $(X,d)$, say, to $\tilde{x}$, since $(X,d)$ is complete, i.e., $\displaystyle \lim_{n\to \infty}f_n(x)=\tilde{x}$
		\\
		Now define $f:X\to X$ as $f(x)=\tilde{x}$.
		Then,
		$$\begin{aligned}
		d(f(x),f(y))&=d(\tilde{x},\tilde{y})\\
		&=d(\lim_{n\to \infty}f_n(x),\lim_{n\to \infty}f_n(y))\\
		&= \lim_{n\to \infty}d(f_n(x),f_n(y))\\
		&\le \lim_{n\to \infty}\Big[c_n \cdot d(x,y)\Big]= \Big[\lim_{n\to \infty}c_n\Big]\cdot d(x,y)\\
		\therefore d(f(x),f(y))&\le c\cdot d(x,y), \quad \text{ where } c\in [0,1)  
		\end{aligned}$$
		Hence $(f_n)_{n\ge1}\to f$ in $(Con(X),\bar{d}_\infty)$.
	\end{proof}
\end{theorem}

Now let us define $\mathcal{D}:\mathcal{S}_X^n\times \mathcal{S}_X^n\to \mathbb{R}$  as $\displaystyle\mathcal{D}(\mathbb{S},\mathbb{T})=\min_{\sigma\in S_n}\sum_{i=1}^{n} \bar{d}_{\infty}(f_i,g_{\sigma(i)})$ where $\mathbb{S}=\{X;f_i\}_{i=1}^n$ and $\mathbb{T}=\{X;g_i\}_{i=1}^n$ are $n$-iterated function systems on $X$.
We note that, if $\mathbb{T}$ is minimally ordered with respect to $\mathbb{S}$, then  $\displaystyle\mathcal{D}(\mathbb{S},\mathbb{T})=\sum_{i=1}^{n} \bar{d}_{\infty}(f_i,g_i)$.

We prove that $\mathcal{D}$ is a metric function.

\noindent\textbf{Claim:} $(\mathcal{S}_X^n,\mathcal{D})$ is a metric space.
\\Let $\mathbb{S}=\{X;f_i\}_{i=1}^n$, $\mathbb{T}=\{X;g_i\}_{i=1}^n$ and $\mathbb{U}=\{X;h_i\}_{i=1}^n$ be $n$-iterated function systems on $X$.
We have \[0\le\mathcal{D}(\mathbb{S},\mathbb{S})=\min_{\sigma\in S_n}\sum_{i=1}^{n} \bar{d}_{\infty}(f_i,f_{\sigma(i)})\le \sum_{i=1}^{n} \bar{d}_{\infty}(f_i,f_i)=0 \]
\[\therefore \mathcal{D}(\mathbb{S},\mathbb{S})=0 \]

Now suppose $\mathcal{D}(\mathbb{S},\mathbb{T})=0$. Then for some $\sigma \in S_n$ we get $\bar{d}_\infty(f_i,g_{\sigma(i)})=0, \forall i=1,2,\cdots,n$. Since $\bar{d}_\infty$ is a metric, we obtain $f_i=g_{\sigma(i)}, \forall i=1,2,\cdots,n.$ Hence, $\mathbb{S}=\mathbb{T}$.

To prove the symmetry property, we have
\begin{align*}\mathcal{D}(\mathbb{S},\mathbb{T})&=\min_{\sigma\in S_n}\sum_{i=1}^{n} \bar{d}_{\infty}(f_i,g_{\sigma(i)})\\&=\min_{\sigma\in S_n}\sum_{i=1}^{n} \bar{d}_{\infty}(g_{\sigma(i)},f_i)
\\
\end{align*}
\begin{align*}
&=\min_{\tau\in S_n}\sum_{i=1}^{n} \bar{d}_{\infty}(g_i,f_{\tau(i)})
\\
&
=\mathcal{D}(\mathbb{T},\mathbb{S})\end{align*}

%
Further to prove $\mathcal{D}(\mathbb{S},\mathbb{U})\le \mathcal{D}(\mathbb{S},\mathbb{T})+\mathcal{D}(\mathbb{T},\mathbb{U})$, i.e., the triangle inequality for $\mathcal{D}$ in $\mathcal{S}_X^n$, without loss of generality we prove the triangle inequality for $\mathcal{D}$ in $\mathcal{S}_X^2$.
\\
For this purpose let $\mathbb{S}=\{X;f_1,f_2\},\mathbb{T}=\{X;g_1,g_2\},\mathbb{U}=\{X;h_1,h_2\}$ and $\mathbb{S},\mathbb{T},\mathbb{U}\in \mathcal{S}_X^2$. Also, let $\mathbb{T}$ with respect to $\mathbb{S}$ and $\mathbb{U}$ with respect to $\mathbb{T}$ be  minimally ordered and $\bar{d}_\infty(f_i,g_j)=\alpha_{ij}, \bar{d}_\infty(g_i,h_j)=\beta_{ij},\bar{d}_\infty(f_i,h_j)=\gamma_{ij}$ for $i,j\in\{1,2\}$ where $\alpha_{ij},\beta_{ij},\gamma_{ij}\in \mathbb{R}$. Then $\alpha_{11}+\alpha_{22}\le \alpha_{12}+\alpha_{21}$ and $\beta_{11}+\beta_{22}\le \beta_{12}+\beta_{21}$. Since, $\bar{d}_\infty$ is a metric, the triangle inequality of $\bar{d}_\infty$ gives $\gamma_{11}\le \alpha_{11}+\beta_{11}$ and $\gamma_{22}\le \alpha_{22}+\beta_{22}$. Therefore, $\gamma_{11}+\gamma_{22}\le (\alpha_{11}+\alpha_{22})+(\beta_{11}+\beta_{22})$.
\\
We have, either $\mathcal{D}(\mathbb{S},\mathbb{U})=\gamma_{11}+\gamma_{22}$ or $\mathcal{D}(\mathbb{S},\mathbb{U})=\gamma_{12}+\gamma_{21}$.
\\
If $\mathcal{D}(\mathbb{S},\mathbb{U})=\gamma_{11}+\gamma_{22}$, then $\mathcal{D}(\mathbb{S},\mathbb{U})\le \mathcal{D}(\mathbb{S},\mathbb{T})+\mathcal{D}(\mathbb{T},\mathbb{U})$, because $\gamma_{11}+\gamma_{22}\le (\alpha_{11}+\alpha_{22})+(\beta_{11}+\beta_{22})$ and $\mathcal{D}(\mathbb{S},\mathbb{T})=\alpha_{11}+\alpha_{22}, \mathcal{D}(\mathbb{T},\mathbb{U})=\beta_{11}+\beta_{22}$.
\\
On the contrary, if $\mathcal{D}(\mathbb{S},\mathbb{U})=\gamma_{12}+\gamma_{21}$, then $\gamma_{12}+\gamma_{21}\le \gamma_{11}+\gamma_{22}\le (\alpha_{11}+\alpha_{22})+(\beta_{11}+\beta_{22})$. Therefore $\mathcal{D}(\mathbb{S},\mathbb{U})\le \mathcal{D}(\mathbb{S},\mathbb{T})+\mathcal{D}(\mathbb{T},\mathbb{U})$.
\\
Hence the triangle inequality for $\mathcal{D}$ holds in $\mathcal{S}_X^2$ and the same follows for $\mathcal{S}_X^n$ by the principle of mathematical induction.
\\
Thus, $\mathcal{D}$ is a metric on $\mathcal{S}_X^n$.
\hfill \qed

In the next example we illustrate the previously discussed concepts.
\begin{example}
	Consider $X=[0,1]$ with the Euclidean metric and $\mathcal{S}_{[0,1]}^2$ with the metric $\mathcal{D}$. Let $\mathbb{S}=\{[0,1];f_1,f_2\},\mathbb{T}=\{[0,1];g_1,g_2\},\mathbb{U}=\{[0,1];h_1,h_2\}$ be iterated function systems where \begin{align*}
	f_1(x)=\frac{1}{2}x,&\quad f_2(x)=\frac{1}{2}+\frac{1}{2}x,\\
	g_1(x)=\frac{1}{3}x,&\quad g_2(x)=\frac{2}{3}+\frac{1}{3}x,\\ h_1(x)=\frac{1}{2}+\frac{1}{2}x,&\quad h_2(x)=\frac{3}{4}x.
	\end{align*} 
	$\begin{aligned}
	\text{We have }\quad \bar{d}_\infty(f_1,g_1)&=\sup_{x\in X}\frac{d(f_1(x),g_1(x))}{1+d(f_1(x),g_1(x))}\\&=\sup_{x\in [0,1]}\frac{|\frac{1}{2}x-\frac{1}{3}x|}{1+|\frac{1}{2}x-\frac{1}{3}x|}=\sup_{x\in [0,1]}\frac{\frac{1}{6}x}{1+\frac{1}{6}x}=\frac{\frac{1}{6}}{1+\frac{1}{6}}=\frac{1}{7},\\
	\end{aligned}$
	$
	\\
	\displaystyle \bar{d}_\infty(f_2,g_2)=\sup_{x\in [0,1]}\frac{|\frac{1}{2}+\frac{1}{2}x-(\frac{2}{3}+\frac{1}{3}x)|}{1+|\frac{1}{2}+\frac{1}{2}x-(\frac{2}{3}+\frac{1}{3}x)|}=\sup_{x\in [0,1]}\frac{|\frac{1}{6}x-\frac{1}{6}|}{1+|\frac{1}{6}x-\frac{1}{6}|}=\frac{\frac{1}{6}}{1+\frac{1}{6}}=\frac{1}{7},
	\\
	\bar{d}_\infty(f_1,g_2)=\sup_{x\in [0,1]}\frac{|\frac{1}{2}x-(\frac{2}{3}+\frac{1}{3}x)|}{1+|\frac{1}{2}x-(\frac{2}{3}+\frac{1}{3}x)|}=\sup_{x\in [0,1]}\frac{|\frac{1}{6}x-\frac{2}{3}|}{1+|\frac{1}{6}x-\frac{2}{3}|}=\frac{\frac{2}{3}}{1+\frac{2}{3}}=\frac{2}{5},
	\\
	\bar{d}_\infty(f_2,g_1)=\sup_{x\in [0,1]}\frac{|\frac{1}{2}+\frac{1}{2}x-\frac{1}{3}x|}{1+|\frac{1}{2}+\frac{1}{2}x-\frac{1}{3}x|}=\sup_{x\in [0,1]}\frac{|\frac{1}{6}x+\frac{1}{2}|}{1+|\frac{1}{6}x+\frac{1}{2}|}=\frac{\frac{2}{3}}{1+\frac{2}{3}}=\frac{2}{5},
	\\
	\bar{d}_\infty(f_1,h_1)=\sup_{x\in [0,1]}\frac{|\frac{1}{2}x-(\frac{1}{2}+\frac{1}{2}x)|}{1+|\frac{1}{2}x-(\frac{1}{2}+\frac{1}{2}x)|}=\frac{\frac{1}{2}}{1+\frac{1}{2}}=\frac{1}{3},	
	\\
	\bar{d}_\infty(f_2,h_2)=\sup_{x\in [0,1]}\frac{|\frac{1}{2}+\frac{1}{2}x-\frac{3}{4}x|}{1+|\frac{1}{2}+\frac{1}{2}x-\frac{3}{4}x|}=\sup_{x\in [0,1]}\frac{|\frac{1}{2}-\frac{1}{4}x|}{1+|\frac{1}{2}-\frac{1}{4}x|}=\frac{\frac{1}{2}}{1+\frac{1}{2}}=\frac{1}{3},
	\\
	\bar{d}_\infty(f_1,h_2)=\sup_{x\in [0,1]}\frac{|\frac{1}{2}x-\frac{3}{4}x|}{1+|\frac{1}{2}x-\frac{3}{4}x|}=\sup_{x\in [0,1]}\frac{\frac{1}{4}x}{1+\frac{1}{4}x}=\frac{\frac{1}{4}}{1+\frac{1}{4}}=\frac{1}{5},
	\\
	\bar{d}_\infty(f_2,h_1)=\sup_{x\in [0,1]}\frac{|\frac{1}{2}+\frac{1}{2}x-(\frac{1}{2}+\frac{1}{2}x)|}{1+|\frac{1}{2}+\frac{1}{2}x-(\frac{1}{2}+\frac{1}{2}x)|}=0.
	$\\
	Also $
	\displaystyle\bar{d}_\infty(f_1,g_1)+\bar{d}_\infty(f_2,g_2)=\frac{1}{7}+\frac{1}{7}=\frac{2}{7},\quad
	\bar{d}_\infty(f_1,g_2)+\bar{d}_\infty(f_2,g_1)=\frac{2}{5}+\frac{2}{5}=\frac{4}{5}
	$ and \linebreak$\displaystyle
	\bar{d}_\infty(f_1,h_1)+\bar{d}_\infty(f_2,h_2)=\frac{1}{3}+\frac{1}{3}=\frac{2}{3},$ $\displaystyle\quad \bar{d}_\infty(f_1,h_2)+\bar{d}_\infty(f_2,h_1)=\frac{1}{5}+0=\frac{1}{5}.
	$
	\\
	Hence $\mathbb{T}$ is minimally ordered with respect to $\mathbb{S}$ and $\mathcal{D}(\mathbb{S},\mathbb{T})=\frac{2}{7}$. But $\mathbb{U}$ is not minimally ordered with respect to $\mathbb{S}$, because $\bar{d}_\infty(f_1,h_2)+\bar{d}_\infty(f_2,h_1)<\bar{d}_\infty(f_1,h_1)+\bar{d}_\infty(f_2,h_2)$. Therefore, the minimal order of $\mathbb{U}$ with respect to $\mathbb{S}$ is $\mathbb{U}=\{[0,1];h_2,h_1\}$ and $\mathcal{D}(\mathbb{S},\mathbb{U})=\frac{1}{5}$.\\
	First few stages of finding the attractors of $\mathbb{S}, \mathbb{T}, \mathbb{U}$ are as shown in Figure \ref{Fig:dist-attr}. Here, while evaluating the attractors of $\mathbb{S}$ and $\mathbb{U}$ we have shifted different iterates a bit to see the overlapping part, if any. Thus, from Figure \ref{Fig:dist-attr} we can see that $\mathbb{S}$ is a just-touching IFS, $\mathbb{T}$ is a totally-disconnected IFS, and $\mathbb{U}$ is an overlapping IFS. Also, note that, even though the attractors of $\mathbb{S}$ and $\mathbb{U}$ are the same, they have a positive metric distance in $\mathcal{D}$.
	
	\begin{figure}[h]
		\centering
		
		\hfill
		\subfigANDtitle{\includegraphics[width=.3\textwidth]{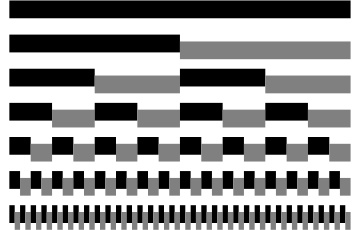} \\ Iterates of $\mathbb{S}$} \hfill
		\subfigANDtitle{\includegraphics[width=.3\textwidth]{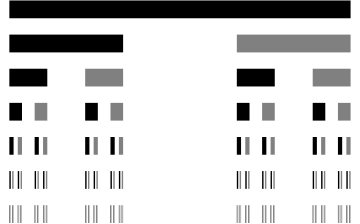} \\ Iterates of $\mathbb{T}$} \hfill\mbox{}
		\\
		\hfill
		\medskip
		\subfigANDtitle{\includegraphics[width=.3\textwidth]{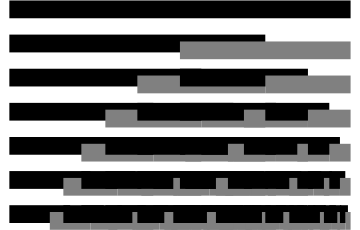} \\ Iterates of $\mathbb{U}$}
		\hfill\mbox{}
		\caption{First six iterates in finding the attractors of $\mathbb{S},\mathbb{T}$ and $\mathbb{U}$}
		\label{Fig:dist-attr}
	\end{figure}
\end{example}
Now we pose the following problem.
Let $\mathbb{T}$ with respect to $\mathbb{S}$ and $\mathbb{U}$ with respect to $\mathbb{T}$ both be minimally ordered. Does it follow that $\mathbb{U}$ with respect to $\mathbb{S}$ is minimally ordered?\\
We found that the implication doesn't follow. The counterexample to this is as follows:
\begin{example}
	Consider the metric space $X=\mathbb{R}^2$ with the Euclidean metric. Let $\mathbb{S}=\{X;s_1,s_2 \},$ $\mathbb{T}=\{X;t_1,t_2 \}$ and $\mathbb{U}=\{X;u_1,u_2 \}$ be iterated function systems on $X$ where $s_1(x,y)=(0,0),$ $s_2(x,y)=(1,0),$ $t_1(x,y)=(0,1),$ $t_2(x,y)=(1,-1),$ $u_1(x,y)=(1,1)$ and $u_2(x,y)=(0,-1)$. The illustration finding the minimal orders are given in Figure \ref{Fig:minimal}.
	\vspace*{.2cm}
	
	We found that $d(s_1,t_1)=1$, $d(s_2,t_2)=1$, $d(s_1,t_2)=\sqrt{2}$ and $d(s_2,t_1)=\sqrt{2}$. Therefore, $\mathbb{T}$ is minimally ordered with respect to $\mathbb{S}$. Also, $d(t_1,u_1)=1$, $d(t_2,u_2)=1$, $d(t_1,u_2)=\sqrt{2}$ and $d(t_2,u_1)=\sqrt{2}$. Hence, $\mathbb{U}$ is minimally ordered with respect to $\mathbb{T}$. But, we get $d(s_1,u_1)=2$, $d(s_2,u_2)=2$, $d(s_1,u_2)=1$ and $d(s_2,u_1)=1$. Thus, $\mathbb{U}$ is minimally ordered with respect to $\mathbb{S}$ and the minimal order of $\mathbb{U}$ with respect to $\mathbb{S}$ is $\mathbb{U}=\{X;u_2,u_1 \}$. Hence, we proved that even though $\mathbb{T}$ with respect to $\mathbb{S}$ and $\mathbb{U}$ with respect to $\mathbb{T}$ both are minimally ordered, $\mathbb{U}$ with respect to $\mathbb{S}$ may not be minimally ordered. Thus, minimally ordered relation is not transitive.
	
	\begin{figure}[h]
		\centering
		\subfigANDtitle{
			\begin{tikzpicture}[scale=.7]
			\node[node2] (1) [label={[label distance=0cm]90:$(0,1)$},label={[label distance=0cm]180:$s_1$}]at (0,3){};
			\node[main node] (3) [label={[label distance=0cm]0:$(1,0)$},label={[label distance=0cm]180:$t_2$}] at (3,0){};
			\node[main node] (4) [label={[label distance=0cm]-180:$(0,0)$},label={[label distance=0cm]0:$t_1$}] at (0,0) {};
			\node[node2] (6) [label={[label distance=0cm]-90:$(1,-1)$},label={[label distance=0cm]0:$s_2$}] at (3,-3) {};
			
			\path[draw,thick]
			(1) edge node [ pos=0.5, sloped, above] {$\sqrt{2}$} (3)
			(1) edge node [pos=0.5, sloped, below]{$1$} (4)
			(3) edge node [ pos=0.5, sloped, above]{$1$}  (6)
			(4) edge node [ pos=0.5, sloped, below]{$\sqrt{2}$} (6);
			
			\end{tikzpicture} \\ (i) $\mathbb{T}$ w.r.t $\mathbb{S}$}
		\hfill\mbox{}
		\hfill
		\subfigANDtitle{
			\begin{tikzpicture}[scale=.7]
			\node[node1] (2) [label={[label distance=0cm]90:$(1,1)$},label={[label distance=0cm]0:$u_2$}] at (3,3){};
			\node[main node] (3) [label={[label distance=0cm]0:$(1,0)$},label={[label distance=0cm]180:$t_2$}] at (3,0){};
			\node[main node] (4) [label={[label distance=0cm]-180:$(0,0)$},label={[label distance=0cm]0:$t_1$}] at (0,0) {};
			\node[node1] (5) [label={[label distance=-0cm]-90:$(0,-1)$},label={[label distance=0cm]180:$u_1$}] at (0,-3) {};
			
			\path[draw]
			(4) edge node [ pos=0.5, sloped, above]{$1$} (5)
			(4) edge node [ pos=0.5, sloped, above]{$\sqrt{2}$} (2)
			(3) edge node [ pos=0.5, sloped, below]{$1$} (2)
			(3) edge node [ pos=0.5, sloped, below]{$\sqrt{2}$} (5);
			
			\end{tikzpicture}\\ (ii) $\mathbb{U}$ w.r.t $\mathbb{T}$}
		\hfill\mbox{}
		\hfill
		\subfigANDtitle{
			\begin{tikzpicture}[scale=.7]
			\node[node2] (1) [label={[label distance=0cm]90:$(0,1)$},label={[label distance=0cm]180:$s_1$}]at (0,3){};
			\node[node1] (2) [label={[label distance=0cm]90:$(1,1)$},label={[label distance=0cm]0:$u_2$}] at (3,3){};
			\node[node1] (5) [label={[label distance=-0cm]-90:$(0,-1)$},label={[label distance=0cm]180:$u_1$}] at (0,-3) {};
			\node[node2] (6) [label={[label distance=0cm]-90:$(1,-1)$},label={[label distance=0cm]0:$s_2$}] at (3,-3) {};
			
			\path[draw,thick]
			(2) edge node [ pos=0.5, sloped, below]{$2$} (6)
			(5) edge node [ pos=0.5, sloped, below]{$1$} (6)
			(2) edge node [ pos=0.5, sloped, above]{$1$} (1)
			(1) edge node [ pos=0.5, sloped, below]{$2$} (5);
			
			\end{tikzpicture} \\ (iii) $\mathbb{U}$ w.r.t $\mathbb{S}$}
		\hfill\mbox{}
		\caption{Finding the minimal order of (i) $\mathbb{T}$ w.r.t $\mathbb{S}$ (ii) $\mathbb{U}$ w.r.t $\mathbb{T}$ and (iii) $\mathbb{U}$ w.r.t $\mathbb{S}$}
		\label{Fig:minimal}
	\end{figure}
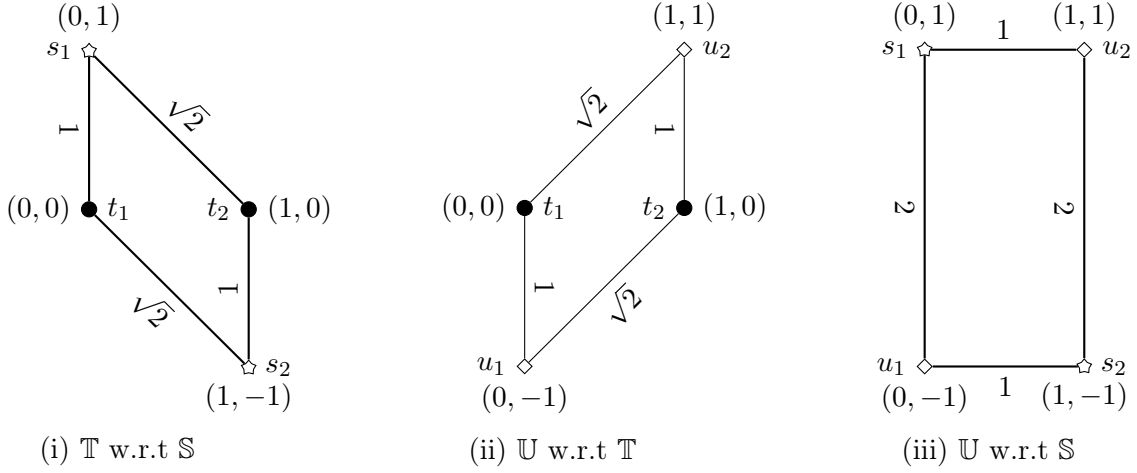
\end{example}
We make the following remark on relation of an $n$-IFS minimally ordered with respect to another $n$-IFS.
\begin{remark}
	The relation of an IFS minimally ordered with respect to another IFS in $\mathcal{S}_X^n$ is reflexive and symmetric, but not transitive.
\end{remark}
The following definition is introduced to discuss on the smaller subsets of $\mathcal{S}_X^n$ in which the relation used in the above remark is an equivalence relation.
\begin{definition}
	A subset $\mathcal{B}$ of $\mathbb{S}_X^n$ is called minimally ordered (abbreviated as $m.\,o$) if the relation of an IFS minimally ordered with respect to another IFS is an equivalence relation on $\mathcal{B}$.
\end{definition}
We remark the following on the $n$- iterated function systems and the minimally ordered sets containing them.
\begin{remark}
	For every $\mathbb{S}\in \mathcal{S}_X^n$, there exists $\mathcal{B}\subset \mathcal{S}_X^n$ such that $\mathcal{B}$ is minimally ordered.
\end{remark}
Now we define the notion of Cauchy and convergent sequences in the setting of sequences of $n$- iterated function systems.
\begin{definition}
	Let $(\mathbb{S}_j)_{j\ge1}$ be a sequence of $n$-iterated function systems on $X$. Then $(\mathbb{S}_j)_{j\ge1}$ is said to be Cauchy if for every $\epsilon>0$, there exists $N\in \mathbb{N}$ such that for all $j,k\ge N$, $\mathcal{D}(\mathbb{S}_j,\mathbb{S}_k)<\epsilon$. Also, we say $(\mathbb{S}_j)_{j\ge1}$ converges to an iterated function system $\mathbb{S}\in \mathcal{S}_X^n$ if for every $\epsilon>0$, there exists $N\in\mathbb{N}$ such that for every $j\ge N$, $\mathcal{D}(\mathbb{S}_j,\mathbb{S})<\epsilon$.
	
\end{definition}
\begin{remark}
	
	A sequence of IFSs $(\mathbb{S}_j)_{j\ge1}$ is said to be convergent if there exists contractions $f_i$ on $X$ such that $(f_{i_j})_{j\ge1}$ converges to $f_i$ point-wise for $i=1,2,\cdots$
\end{remark}

%

\begin{theorem}\label{conv3}
	Consider the $n$-iterated function systems $\mathbb{S}_j\in \mathcal{S}_X^n$ for $j=1,2,\cdots$. Assume $(\mathbb{S}_j)_{j\ge1}\to \mathbb{S}$ in $(\mathcal{S}_X^n,\mathcal{D})$. Let $A_j$ be the attractor of the iterated function systems $\mathbb{S}_j$ for $j=1,2,\cdots$ and $A$ be the attractor of the iterated function system $\mathbb{S}$. Then $A_j\to A$ in the Hausdorff metric on $X$.
	\begin{proof}
		We prove $f_{i_j}\to f_{\sigma_0(i)}$ in the metric $\bar{d}$ for some $\sigma_0\in S_n$.
		Let $\epsilon>0$ be given. Since $\mathbb{S}_j\to \mathbb{S}$, there exists $N\in \mathbb{N}$ such that for all $j\ge N$, $\mathcal{D}(\mathbb{S}_j,\mathbb{S})<\epsilon$. Then $\displaystyle\min_{\sigma\in S_n}\sum_{i=1}^n\bar{d}_\infty(f_{i_j},f_{\sigma(i)})<\epsilon$.
		Then $\displaystyle\sum_{i=1}^n\bar{d}_\infty(f_{i_j},f_{\sigma_0(i)})<\epsilon$,
		where $\sigma_0\in S_n$ is such that $\displaystyle\sum_{i=1}^n\bar{d}_\infty(f_{i_j},f_{\sigma(i)})$ is minimum for $\sigma=\sigma_0$. Hence we have $\bar{d}_\infty(f_{i_j},f_{\sigma_0(i)})<\epsilon$ for every $j\ge N$. Therefore $f_{i_j}\to f_{\sigma_0(i)}$ in the metric $\bar{d}_\infty$. 
	\end{proof}
\end{theorem} 
The following is an alternate proof for the Theorem \ref{conv3}.
\begin{proof}
	Let $\Lambda$ be the code space on $\{1,2,\cdots,n\}$ 
	and $\gamma=\gamma_1\gamma_2\gamma_3\cdots\in \Lambda$. 
	Suppose $\pi_j$ and $\pi$ be the code maps for $A_{j}$ and $A$ respectively.\\
	Then $\pi_j(\gamma):=f_{\gamma_j}(x):=\cdots\circ f_{\gamma_{3_j}}\circ f_{\gamma_{2_j}} \circ f_{\gamma_{1_j}}(x)$ for every $j\in\mathbb{N}$ and $\displaystyle\pi(\gamma):=f_\gamma(x):=\lim_{n\to \infty}f_{\gamma_n}\circ\cdots\circ f_{\gamma_3}\circ f_{\gamma_2}\circ f_{\gamma_1}(x)$ for any $x\in X$.\\
	We have $f_{\gamma_{k_j}}\to f_{\gamma_k}$. Therefore $f_{\gamma_j}\to f_\gamma$ and $f_{\gamma_j}(x)\to f_\gamma(x)$ for every $x\in X$ and $\gamma\in \Lambda$.\\
	But we have $\displaystyle A_{j}=\bigcup_{x\in X,\gamma\in\Lambda}f_{\gamma_j}(x)$ and $\displaystyle A=\bigcup_{x\in X,\gamma\in\Lambda}f_{\gamma}(x)$.\\
	Hence $A_{j}\to A$ in the Hausdorff metric on $X$.
\end{proof}
The following  example is a demonstration of the theorem \ref{conv3}.
\begin{example}
	Let $\mathbb{S}_j=\{[0,1]; \frac{1}{3}x+\frac{1}{3j}, \frac{2}{3}+\frac{1}{3}x\}$ and $\mathbb{S}=\{X;\frac{1}{3}x,\frac{2}{3}+\frac{1}{3}x\}$.
	Then $(\mathbb{S}_j)_{j\ge1}\to \mathbb{S}$. A demonstration of $A_{\mathbb{S}_j}\to A_{\mathbb{S}}$ is provided in the Figure \ref{Fig:attr}. We note that, in $(\mathbb{S}_n)_{n\ge1}$ only $\mathbb{S}_1$ is just touching and $\mathbb{S}_n$ are totally disconnected for all $n\ge 2$.
	\begin{figure}[h]
		\centering
		
		\hfill
		\subfigANDtitle{\includegraphics[width=.3\textwidth]{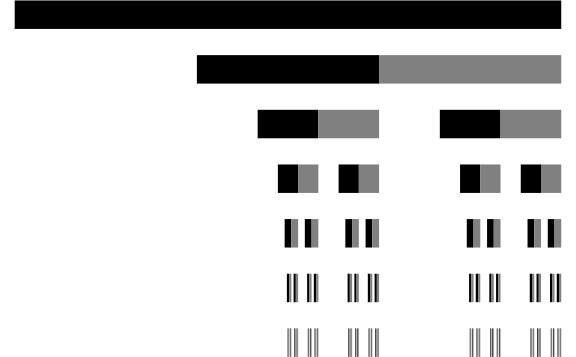} \\ Iterates of $\mathbb{S}_1$} \hfill
		\subfigANDtitle{\includegraphics[width=.3\textwidth]{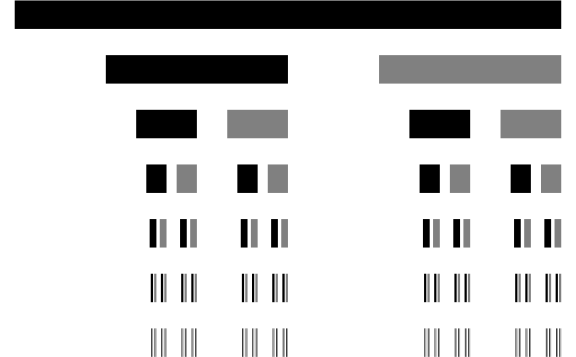} \\ Iterates of $\mathbb{S}_2$}
		\hfill\mbox{}
		
		\medskip
		
		\hfill
		\subfigANDtitle{\includegraphics[width=.3\textwidth]{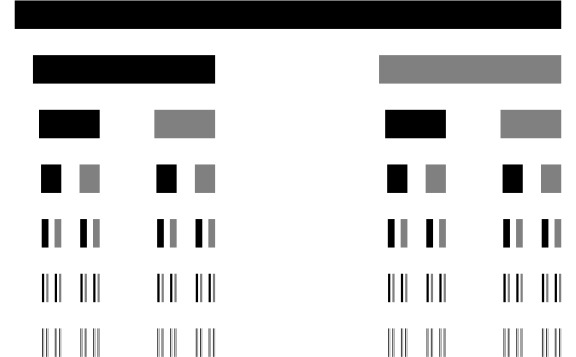} \\ Iterates of  $\mathbb{S}_{10}$} \hfill
		\subfigANDtitle{\includegraphics[width=.3\textwidth]{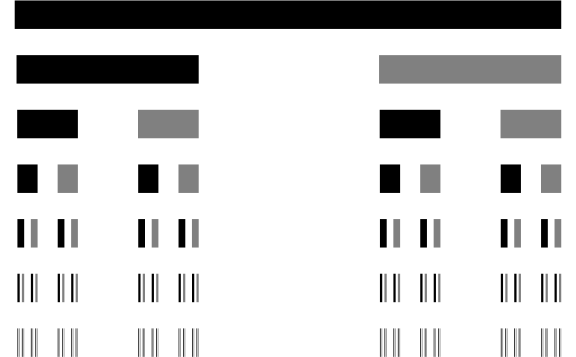} \\ Iterates of $\mathbb{S}_{100}$}
		\hfill\mbox{}
		\medskip
		
		\hfill
		\subfigANDtitle{\includegraphics[width=.3\textwidth]{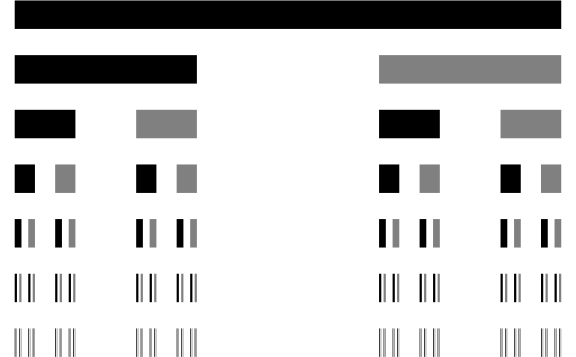} \\ Iterates of $\mathbb{S}$}
		\hfill\mbox{}
		\caption{First few iterates in finding the attractors of $(\mathbb{S}_j)_{j\ge1}$ and $\mathbb{S}$}
		\label{Fig:attr}
	\end{figure}
\end{example}

\begin{proposition}
	Consider a complete metric space $(X,d)$. Let $C(X)$ denotes the set of all continuous functions on $X$. Then the metric space $(C(X),\bar{d}_{\infty})$ is complete, where $\bar{d}_{\infty}$ is defined by $\displaystyle\bar{d}_{\infty}(f,g)=\sup_{x\in X}\frac{d(f(x),g(x))}{1+d(f(x),g(x))}$ for every $f,g\in C(X)$.
	\begin{proof}
		Let $(f_n)_{n\ge1}$ be a Cauchy sequence in $C(X)$. Consider the sequence $(f_n(x))_{n\ge1}$ for some $x\in X$. We will prove $(f_n(x))_{n\ge1}$ is a Cauchy sequence. Let $\epsilon>0$ be given. Then, $(f_n)$ being a Cauchy sequence, there exists $N\in\mathbb{N}$ such that $\bar{d}_\infty(f_n,f_m)<\epsilon $, for every $n,m\ge N$. Thus, $\displaystyle\frac{d(f_n(x),f_m(x))}{1+d(f_n(x),f_m(x))}<\epsilon$, for every $n,m\ge N$ and $x\in X$ (here $N$ does not depend on $x$). Hence $d(f_n(x),f_m(x))<\frac{\epsilon}{1-\epsilon}$. Therefore, $f_n(x)$ is a Cauchy sequence in $(X,d)$. But $(X,d)$ is complete, and hence $f_n(x)$ converges in $(X,d)$, to say $x'\in X$. Now define $f:X\to X$ by $f(x)=x'$ for every $x\in X$. Then $f_n(x)\to f(x)$ in $(X,d)$ for each $x\in X$. Thus $f_n\to f$ uniformly.
		
		Now for a given $\epsilon>0$, there exists $\delta_n>0$ for each $n\in \mathbb{N}$ such that, \linebreak$d(f_n(x),f_n(y))<\epsilon$  whenever $d(x,y)<\delta_n$. Then, on taking $\displaystyle\delta = \inf_{n\in\mathbb{N}}\delta_n$ (infimum exists since $(\delta_n)_{n\ge1}$ is bounded below by $0$), we get, if $d(x,y)<\delta$ then $$\begin{aligned}d(f(x),f(y))&=d(\lim_{n\to \infty}f_n(x),\lim_{n\to \infty}f_n(y))\\
		&=\lim_{n\to \infty}d(f_n(x),f_n(y)) \quad [\because d:X\times X\to \mathbb{R}\text{ is continuous}]\\
		&\le \epsilon
		\end{aligned}$$
		Hence $f\in C(X)$.
		Thus, if $(X,d)$ is complete, then $(C(X),\bar{d}_\infty)$ is also complete.
		
	\end{proof}
\end{proposition}

%

\begin{theorem}
	A minimally ordered eventually decreasing, Cauchy sequence of hyperbolic iterated function systems with $n$-contractions on a complete metric space $(X,d)$ is convergent in $(\mathcal{S}_X^n,\mathcal{D})$.
	\begin{proof}
		Let $(\mathbb{S}_j)_{j\ge1}$ be an eventually decreasing Cauchy sequence of $n$-iterated function systems on a complete metric space $(X,d)$, where $\mathbb{S}_j=\{X;f_{i_j}\}_{i=1}^{n}$ for $j=1,2,\cdots$ Then $(f_{i_j})_{j\ge1}$ is an eventually decreasing Cauchy sequence of contractions on the complete metric space $(X,d)$. Then by Theorem \ref{con}, the sequence converges, say, to $f_i\in Con(X)$.\\ Now, let $\mathbb{S}=\{X;f_i\}_{i=1}^n$. Then $\mathbb{S}\in \mathcal{S}_X^n$.
		\\
		\textbf{Claim:} $(\mathbb{S}_j)_{j\ge1}\to \mathbb{S}$ in $(\mathcal{S}_X^n,\mathcal{D})$.
		
		Let $\epsilon>0$ be given. Then, since $(f_{i_j})_{j\ge1}\to f_i$, for each $i=1,2,\cdots,n$, there exists $N_i\in\mathbb{N}$ such that $\displaystyle\bar{d}_\infty(f_{i_j},f_i)<\frac{\epsilon}{n}$, for every $j\ge N_i$. Take $\displaystyle N=\max_{i=1}^{n}{N_i}$.\\
		Then ,
		$\displaystyle\mathcal{D}(\mathbb{S}_j,\mathbb{S})=\sum_{i=1}^n\bar{d}_\infty(f_{i_j},f_i)<n\cdot \frac{\epsilon}{n}=\epsilon$ for every $j\ge N$. \\
		Therefore, $(\mathbb{S}_j)_{j\ge1}\to \mathbb{S}$ in $(\mathcal{S}_X^n,\mathcal{D})$.
		
	\end{proof}
\end{theorem}
Some potential applications of the theory developed in this article are discussed in the following section.
\section{Applications}
In this section we propose certain ways to apply the theory of sequences of iterated function systems developed. Suppose we have the data images of the geographical location of a soon to be extinct species in a world map at varying times. Assume that the time gap between taking images is decreasing, and the corresponding series of time gaps is convergent, say to time $T$. With the given images, using Collage theorem we construct a finite sequence of iterated function systems. We find a decreasing sequence of iterated function system whose speed of convergence matches with the speed convergence of the sequence of time gaps in taking the images, which can be controlled. The attractor of the IFS to which the sequence of iterated function systems converges will have an approximation of the image at the time $T_0+T$, where $T_0$ is the first instance when we record an image. Hence we can obtain an approximation of the image of geographical data of the species at a later time. Now, suppose we take some measures to increase the population of this species at various locations and record further images and incorporate it to our previous process, then at time $T_0+T$ we will get an attractor showing the effect of measures taken to increase the population of the species. It could be possible to regulate the measures to get the final image at $T_0+T$ as we desire.
The same procedure may be applied to get the effect of medicines on the human body.

We will discuss now a practical scenario in which we can make use of the theory in a much simpler way. According to the $2008$ global assessment conducted by International Union of Conservation of Nature(IUCN), the asian elephant(\textit{Elephas Maximus}) is in the red list and is marked as endangered. The geographical range of the species according to the IUCN sources is as in the \textit{Figure}. To apply our theory in this scenario, we first have to find an $IFS$, say $IFS_1$, using the Collage theorem of fractal theory so that the attractor of $IFS_1$ is $\epsilon$ close to the original image of geographical ranges in the Hausdorff metric. Suppose we are able to find the geographical range in the subsequent years till the current year. Following the same procedure we can find $IFS_2,$ $IFS_3,\cdots ,$ $IFS_{12}$ for the geographical range of the years $2009,\,2010,\cdots,\,2019$. Now we find the convergence pattern for the functions in the $IFSs$: $IFS_1,\,IFS_2,\cdots,\,IFS_{12}$ and manipulate a convergent sequence of iterated function systems $(IFS_n)_{n\ge1}$ with fixed number of contractions. We can use this sequence to find the approximate image of geographical ranges in the future years, simply by finding the attractor of the term in the sequence corresponding to the year we need. For example, in order to find an approximation to the geographical range for the year $2025$, we just need to find the attractor of the $IFS_{18}$ from the sequence $(IFS_n)_{n\ge1}$.
\begin{figure}
	\begin{center}
	\includegraphics[width=15cm]{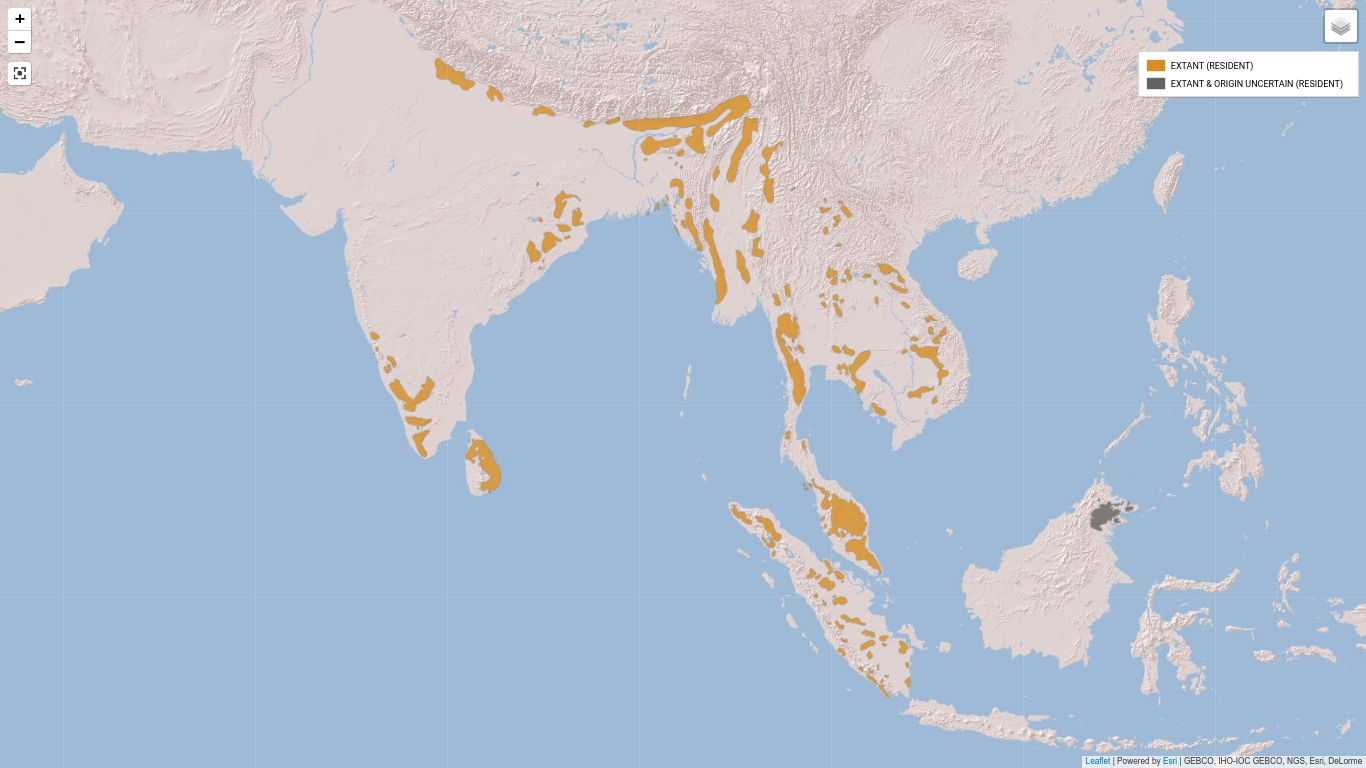}
\caption{Geographical range of Elephas maximus measured in the year $2008$\cite{IUCN:2008}}
\end{center}
\end{figure}

We conclude the findings of our paper in the next section.
\section{Conclusion}
In this paper, we have introduced a sequence of iterated function systems consisting of a fixed number of contractions. In order to achieve this, a proper metric $\mathcal{D}$ is defined on the set of all $n$-iterated function systems on a complete metric space $X$. Further, minimally ordered, decreasing, eventually decreasing, Cauchy type of sequences are discussed, and the convergence of certain types of sequences of $n$-iterated function systems are studied. Towards the end of the article, we have discussed certain applications of the theory developed.
\section*{Acknowledgement}
The first author is very grateful to Council of Scientific \& Industrial Research(CSIR), India for their financial support.
\bibliographystyle{unsrt}
\bibliography{reference}{}
\end{document}